\documentclass[11pt]{article}
\usepackage{amssymb}
\usepackage{amsthm}
\usepackage{multirow}
\usepackage{graphicx}
\usepackage{amsmath}
\usepackage{amsfonts}
\usepackage{mathrsfs}

\newtheorem{Theorem}{Theorem}

\newtheorem{Proposition}{Proposition}

\newtheorem{Corollary}{Corollary}

\title{On the geometry of small weight codewords\\ of dual algebraic geometric codes}
\author{Claudio Fontanari and Chiara Marcolla}
\date{}

\begin{document}
\maketitle

\begin{abstract}
We investigate the geometry of the support of small weight codewords of dual algebraic
geometric codes on smooth complete intersections by applying the powerful tools recently 
developed by Alain Couvreur. In particular, by restricting ourselves to the case of Hermitian codes, we recover and extend previous results obtained by the second named 
author joint with Marco Pellegrini and Massimiliano Sala. 
\end{abstract}

\section{Introduction}

In the recent contribution \cite{Pellegrini}, the number of small weight codewords 
for some families of Hermitian codes is determined. Besides explicit computation, 
the main ingredient in \cite{Pellegrini} is a nice geometric characterization of 
the points in the support of a minimum weight codeword, which turn out to be 
collinear (see Corollary 1 and Proposition 2 in \cite{Pellegrini}). 

Here we show that such a property is not peculiar to Hermitian codes, but it holds 
in full generality for dual algebraic geometric codes on any smooth complete 
intersection projective variety of arbitrary dimension. Namely, by exploiting 
the algebro-geometric tools provided by \cite{C}, we prove the following:

\begin{Theorem}\label{general}
Let $X \subset \mathbb{P}^r$, $r \ge 2$, be a smooth connected complete intersection 
defined over $\mathbb{F}_q$. Let $G$ be a divisor on $X$ such that $L(G) \supseteq 
H^0(\mathbb{P}^r, \mathcal{O}_{\mathbb{P}^r}(m))$. Let $P_1, \ldots, P_n$ be distinct 
rational points on $X$ avoiding the support of $G$ and let $D := P_1 + \ldots + P_n$.   
Let $d$ be the minimum distance of the code $C(D,G)^*$ and let $\{P_{i_1}, \ldots, P_{i_d} \}$ 
be the points in the support of a minimum weight codeword.
\begin{itemize}
\item[(i)] If $d \le m+2$, then $d = m+2$ and all the $m+2$ points $P_{i_j}$ 
are collinear in $\mathbb{P}^r$.
\item[(ii)] If $d \le 2m+2$ and no $m+2$ of the $P_{i_j}$'s are collinear, then $d = 2m+2$ and all the $2m+2$ points 
$P_{i_j}$ lie on a plane conic.
\item[(iii)] If $d \le 3m$, no $m+2$ of the $P_{i_j}$'s are collinear and no $2m+2$ of them lie on a plane conic, then 
$d = 3m$ and all the $3m$ points $P_{i_j}$ lie at the intersection of two complanar plane curves of respective 
degrees $3$ and $m$. 
\end{itemize}
\end{Theorem}

If we focus on the explicit interesting example of Hermitian codes (as presented 
for instance in \cite{Pellik}), then our general result specializes as follows:

\begin{Corollary}\label{distance}
Let $X \subset \mathbb{P}^2$ be the Hermitian curve of affine equation $x^{q+1}=y^q+y$ 
defined over $\mathbb{F}_{q^2}$. Let $G = P_0^\rho$, where $P_0$ is the point at infinity, 
$D=P_1+\ldots+P_n$ with $n=q^3$, and let $C(D,G)^*$ be the Hermitian code. Let $d$ be the 
minimum distance and let $\{P_{i_1}, \ldots, P_{i_d} \}$ be the points in the support of 
a minimum weight codeword.

\begin{itemize}

\item If $0 \leq \rho \leq q^2-q-2$, then all the points $P_{i_j}$ are collinear. 

\item If $\rho \ge q^2-q-2$, then the following holds.
\begin{itemize}
\item[(i)] If $d \le q$, then all the points $P_{i_j}$ are collinear. 
\item[(ii)] If $d \le 2q-2$, then either $q$ of the points $P_{i_j}$ are collinear, or all the points $P_{i_j}$ lie on a plane conic.
\item[(iii)] If $d \le 3q-6$, then either $q$ of the points $P_{i_j}$ 
are collinear, or $2q-2$ of the points $P_{i_j}$ lie on a plane conic, or all the points 
$P_{i_j}$ lie at the intersection of two complanar plane curves of respective degrees $3$ and $q-2$. 
\end{itemize}

\item If $\rho \ge q^2-1$, then the following holds.
\begin{itemize}
\item[(i)] If $d \le q+1$, then all the points $P_{i_j}$ are collinear. 
\item[(ii)] If $d \le 2q$, then either $q+1$ of the points $P_{i_j}$ are collinear, or all the points $P_{i_j}$ lie on a plane conic.
\item[(iii)] If $d \le 3q-3$, then either $q+1$ of the points $P_{i_j}$ 
are collinear, or $2q$ of the points $P_{i_j}$ lie on a plane conic, or all the points 
$P_{i_j}$ lie at the intersection of two complanar plane curves of respective degrees $3$ and $q-1$. 
\end{itemize}

\end{itemize}

\end{Corollary}

In the case of Hermitian codes we may even describe the geometry of small, 
even not minimum, weight codewords (notice that our technical assumption 
$L(G) = H^0(\mathbb{P}^2, \mathcal{O}_{\mathbb{P}^2}(d-2))$ is satisfied 
by the so-called corner codes according to the terminology of \cite{Pellegrini}, 
Definition 2): 

\begin{Proposition}\label{almostminimalweight}
Let $X \subset \mathbb{P}^2$ be the Hermitian curve of affine equation $x^{q+1}=y^q+y$ 
defined over $\mathbb{F}_{q^2}$. Let $G = P_0^\rho$, where $P_0$ is the point at infinity, 
$D=P_1+\ldots+P_n$ with $n=q^3$, and let $C(D,G)^*$ be the Hermitian code. If 
$0 \leq \rho \leq q^2-q-2$ and $L(G) = H^0(\mathbb{P}^2, \mathcal{O}_{\mathbb{P}^2}(d-2))$, then at least $d-1$ of the points $\{P_{i_1}, \ldots, P_{i_{d+a}} \}$ in the support of a codeword of weight $d+a$, where $0 \le a \le d-3$, are collinear. 
\end{Proposition}

In the next future, the second named author is going to apply the above geometric characterization to compute the number of small weight codewords for Hermitian 
corner codes, in the spirit of \cite{Pellegrini}. Indeed, the present research 
is part of her PhD program, supervised by Massimiliano Sala. Both authors are grateful 
to him for strongly stimulating this joint project. The first named author has been 
partially supported by GNSAGA of INdAM and MIUR Cofin 2008 - Geo\-metria 
delle variet\`{a} algebriche e dei loro spazi di moduli (Italy).

\section{The proofs}

\emph{Proof of Theorem \ref{general}.} 
Let $c \in C(D,G)^*$ be a minimum weight codeword having support $\{P_{i_1}, \ldots, P_{i_s} \}$
with $s=d$. By the definition of a dual code, we have 
$$
\sum_{j=1}^s c_j f(P_{i_j})=0
$$
for every $f \in L(G) \supseteq H^0(\mathbb{P}^r, \mathcal{O}_{\mathbb{P}^r}(m))$. 
In particular, we have
$$
\sum_{j=1}^s c_j \mathrm{ev}_{P_{i_j}}(f)=0
$$
for every $f \in H^0(\mathbb{P}^r, \mathcal{O}_{\mathbb{P}^r}(m))$. Hence $\mathrm{ev}_{P_{i_j}}$ 
turn out to be linearly dependent in $H^0(\mathbb{P}^r, \mathcal{O}_{\mathbb{P}^r}(m))^*$ and, 
in the terminology of \cite{C}, Definition 2.4 and Definition 2.8, $P_{i_1}, \ldots, P_{i_s}$
are $m$-linked. Now, we have $d=s=m+2$ in case (i) by \cite{C}, Proposition 6.1,  
$d=s=2m+2$ in case (ii) by \cite{C}, Proposition 7.2, and $d=s=3m$ in case (iii) 
by \cite{C}, Proposition 8.1. Finally, we complete case (i) by \cite{C}, Proposition 6.2, 
case (ii) by \cite{C}, Proposition 7.4, and case (iii) by \cite{C}, Proposition 9.1.

\qed

\noindent
\emph{Proof of Corollary \ref{distance}.} Since $L(G)= \langle \{ x^iy^j: 
i \ge 0, 0 \le j \le q-1, iq+j(q+1) \le \rho \} \rangle$ and $H^0(\mathbb{P}^r, 
\mathcal{O}_{\mathbb{P}^r}(k))$ can be identified with the set of polynomials 
in $r$ variables of degree at most $k$, it is easy to check that 
$L(G) \supseteq H^0(\mathbb{P}^2, \mathcal{O}_{\mathbb{P}^2}(k))$
for every $k \ge 0$ such that $k \le q-1$ and $k(q+1) \le \rho$.\\
If $0 \leq \rho \leq q^2-q-2$, then let $m := d-2$. 
By \cite{Pellik}, \S 5.3, we have $\rho=2q^2-q-uq-v-1$
with $1\leq u, v\leq q-1$ and $d= (q-u)q-v$ if $u < v$, $d=(q-u)q$ if $u\geq v$. 
Hence $d=m+2$ implies $m = u-1$ for $u > v$ and $m = u$ for $u=v$, hence 
\begin{equation}\label{enough}
m(q+1) \le \rho.
\end{equation}

If $\rho \ge q^2-q-2$ let $m := q-2$, if instead $\rho \ge q^2-1$ let 
$m := q-1$. In both cases (\ref{enough}) is satisfied, 
implying $L(G) \supseteq H^0(\mathbb{P}^2, \mathcal{O}_{\mathbb{P}^2}(m))$. 
Now our claim follows from Theorem \ref{general}.

\qed

\noindent
\emph{Proof of Proposition \ref{almostminimalweight}.}
If $c_j$ are the non-zero components of the corresponding codeword, then 
$$
\sum_{j=1}^{d+a} c_j f(P_{i_j})=0
$$
for every $f \in L(G) \supseteq H^0(\mathbb{P}^2, \mathcal{O}_{\mathbb{P}^2}(d-2))$, 
in particular $\{P_{i_1}, \ldots, P_{i_{d+a}} \}$ are $(d-2)$-linked. 

If they are not minimally $(d-2)$-linked, then (up to reordering) we have 
$$
\sum_{j=1}^{d+a-1} b_j f(P_{i_j})=0
$$
for every $f \in H^0(\mathbb{P}^2, \mathcal{O}_{\mathbb{P}^2}(d-2))$. Our assumption
$L(G) = H^0(\mathbb{P}^2, \mathcal{O}_{\mathbb{P}^2}(d-2))$ implies that the $b_j$'s 
are the components of a codeword of weight strictly less than $d+a$. 
By induction on $a$ starting from Corollary \ref{distance} 
we conclude that at least $d-1$ of them are collinear. 

Assume now that the points $\{P_{i_1}, \ldots, P_{i_{d+a}} \}$ are minimally $(d-2)$-linked.
If they are not collinear, there exists a hyperplane $H$ containing exactly $l$ of them, 
with $2 \le l \le d+a-1$. By \cite{C}, Lemma 5.1, the remaining $d+a-l$ points are 
(minimally) $(d-3)$-linked and from \cite{C}, Proposition 7.2, it follows that at least 
$(d-1)$ of them are collinear since by our numerical assumption on $a$ we have 
$d+a-l \le 2(d-3)+1$.

\qed

\hspace{0.5cm}

\noindent
Claudio Fontanari and Chiara Marcolla\newline
Dipartimento di Matematica, Universit\`a di Trento \newline 
Via Sommarive 14, 38123 Trento, Italy. \newline
E-mail addresses: fontanar@science.unitn.it, \newline
chiara.marcolla@unitn.it

\begin{thebibliography}{99}

\bibitem{C} A. Couvreur: The dual minimum distance of arbitrary dimensional algebraic-geometric codes. arXiv:0905.2345 (2009). 

\bibitem{Pellik} T. H{\o}holdt, J. H. van Lint, and R. Pellikaan: 
Algebraic geometry of codes. Handbook of Coding Theory, 1998, 871--961.

\bibitem{Pellegrini} M. Pellegrini, C. Marcolla, and M. Sala: On the weight 
of affine-variety codes and some Hermitian codes. Proceeedings of WCC 2011. 

\end{thebibliography}
\end{document}